\documentclass[english]{amsart}
\usepackage[T1]{fontenc}
\usepackage[latin1]{inputenc}
\usepackage{amssymb}

\makeatletter
 \theoremstyle{plain}    
 \newtheorem{thm}{Theorem}[section]
 \theoremstyle{plain}
 \theoremstyle{definition}
  \newtheorem{example}[thm]{Example}
 \theoremstyle{remark}
 \newtheorem{rem}[thm]{Remark}
 \theoremstyle{plain}    
 \newtheorem{lem}[thm]{Lemma} %%Delete [thm] to re-start numbering
 \theoremstyle{definition}
 \newtheorem{defn}[thm]{Definition}
 \theoremstyle{plain}    
 \newtheorem{prop}[thm]{Proposition} %%Delete [thm] to re-start numbering

\usepackage[all]{xy}

\usepackage{babel}
\makeatother
\begin{document}

\title{Dendroidal Sets}

\author{I. Moerdijk, I. Weiss}

\begin{abstract}
We introduce the concept of a dendroidal set. This is a generalization
of the notion of a simplicial set, specially suited to the study of
operads in the context of homotopy theory. We define a category of
trees, which extends the category $\Delta$ used in simplicial sets,
whose presheaf category is the category of dendroidal sets. We show
that there is a closed monoidal structure on dendroidal sets which
is closely related to the Boardman-Vogt tensor product of operads.
Furthermore we show that each operad in a suitable model category
has a coherent homotopy nerve which is a dendroidal set, extending
another construction of Boardman and Vogt. There is also a notion
of an inner Kan dendroidal set which is closely related to simplicial Kan
complexes. Finally, we briefly indicate the theory of dendroidal objects
and outline several of the applications and further theory of dendroidal
sets. 
\end{abstract}
\maketitle

\section{Introduction}

There is an intimate relation between simplicial sets and categories
(and, more generally, between simplicial objects and enriched categories),
which plays a fundamental role in many parts of homotopy theory. The
goal of this paper is to introduce an extension of the category of
simplicial sets, suitable for studying operads. We call the objects
of this larger category ''dendroidal sets'', and denote the inclusion
functor by \[
i_{!}:(\textrm{simplicial sets})\rightarrow(\textrm{dendroidal sets}).\]

The pair of adjoint functors\[
\xymatrix{\tau:(\textrm{simplicial sets})\ar@<2pt>[r] & (\textrm{categories}):N\ar@<2pt>[l]}
\]
 where $N$ denotes the nerve and $\tau$ its left adjoint, will be
seen to extend to a pair\[
\xymatrix{\tau_{d}:(\textrm{dendroidal sets})\ar@<2pt>[r] & (\textrm{operads}):N_{d}\ar@<2pt>[l]}
\]
having similar properties.

Many other properties and constructions of simplicial sets also extend
to dendroidal sets. In particular, we will show that the cartesian
closed monoidal structure on simplicial sets extends to a (non-cartesian!)
closed monoidal structure on dendroidal sets. Here ''extends'' means
that there is a natural isomorphism\[
i_{!}(X\times Y)\cong i_{!}(X)\otimes i_{!}(Y)\]
for any two simplicial sets $X$ and $Y$. This tensor product of
dendroidal sets is closely related to the Boardman-Vogt tensor product
of operads. In fact, the latter can be defined in terms of the former
by the isomorphism\[
\mathcal{P}\otimes_{BV}\mathcal{Q}\cong\tau_{d}(N_{d}\mathcal{P}\otimes N_{d}\mathcal{Q})\]
for any two operads $\mathcal{P}$ and $\mathcal{Q}$.

We will also define a notion of inner (or weak) Kan complex for dendroidal
sets, extending the simplicial one in the sense that for any simplicial
set $X$, one has that $X$ is an inner Kan complex iff $i_{!}(X)$
is. The nerve of an operad always satisfies this dendroidal inner Kan
condition, just like the nerve of a category satisfies the simplicial
inner Kan condition. Moreover, this inner Kan condition has various
basic properties related to the monoidal structure on dendroidal sets,
the most significant one being that, under some conditions on a dendroidal set $X$,
$\underline{Hom}(X,K)$ is an inner Kan complex whenever $K$ is. The analogous property for simplicial
sets was recently proved by Joyal, and forms one of the basic steps
in the proof of the existence of the closed model structure on simplicial
sets in which the inner Kan complexes are exactly the fibrant objects.
Joyal calls these inner Kan complexes quasi-categories, and one might
call a dendroidal set a \emph{quasi-operad} if it satisfies our dendroidal
version of the inner Kan condition. We expect that there is a closed
model structure on dendroidal sets in which the quasi-operads are
the fibrant objects. Dendroidal sets also seem to be useful in the
theory of homotopy-$\mathcal{P}$-algebras for an operad $\mathcal{P}$
and weak maps between such algebras. We will comment on this towards
the end of this paper. 

The results in this paper were first presented at the Mac Lane Memorial conference in Chicago (April 2006). 
We would like to thank C. Berger, J. Guti\'errez, A. Joyal, A. Lukacs, and M. Shulman for useful comments on early versions.

\section{Operads}

In this paper, \emph{operad} means \emph{coloured symmetric operad}. (In the literature such operads are also referred to as symmetric
multi-categories \cite{Leinster book}.) We briefly
recall the basic definitions, and refer to \cite{BM col operads}
for a more extensive discussion. An operad $\mathcal{P}$ is given
by a set of colours $C$, and for each $n\ge0$ and each sequence
of colours $c_{1},\cdots,c_{n},c$ a set $\mathcal{P}(c_{1},\cdots,c_{n};c)$
(to be thought of as operations taking $n$ inputs of colours $c_{1},\cdots,c_{n}$
respectively to an output of colour $c$). Moreover, there are structure
maps for units and composition. If we write $I=\{*\}$ for the one-point
set, there is for each colour $c$ a unit map \[
u:I\rightarrow\mathcal{P}(c;c)\]
 taking $*$ to $1_{c}$. The composition operations are maps \[
\mathcal{P}(c_{1},\cdots,c_{n};c)\times\mathcal{P}(d_{1}^{1},\cdots,d_{k_{1}}^{1};c_{1})\times\cdots\times\mathcal{P}(d_{1}^{n},\cdots,d_{k_{n}}^{n};c_{n})\rightarrow\mathcal{P}(d_{1}^{1},\cdots,d_{k_{n}}^{n};c)\]
 which we denote $p,q_{1},\cdots,q_{n}\mapsto p(q_{1},\cdots,q_{n})$.
These operations should satisfy the usual associativity and unitary
conditions. Furthermore for each $\sigma\in\Sigma_{n}$ and colours
$c_{1},\cdots,c_{n},c\in C$ there is a map $\sigma^{*}:\mathcal{P}(c_{1},\cdots,c_{n};c)\rightarrow\mathcal{P}(c_{\sigma(1)},\cdots,c_{\sigma(n)};c)$.
These maps define a right action of $\Sigma_{n}$ in the sense that
$(\sigma\tau)^{*}=\tau^{*}\sigma^{*}$, and the composition operations
should be equivariant in some natural sense. The definition can equivalently
be cast in terms of the units and the ''$\circ_{i}$-operations''
\[
\xymatrix{\mathcal{P}(c_{1},\cdots,c_{n};c)\times\mathcal{P}(d_{1},\cdots,d_{k};c_{i})\ar[r]^{{\circ_{i}\,\,\,\,\,\,\,\,\,\,\,}} & \mathcal{P}(c_{1},\cdots,c_{i-1},d_{1},\cdots,d_{k},c_{i+1},\cdots,c_{n};c).}
\]
 A coloured operad $\mathcal{P}$ with set $C$ of colours will also
be referred to as an operad \emph{coloured by} $C$, or an operad
\emph{over} $C$.

The same definition of operad still makes sense if we replace $Set$
by an arbitrary cocomplete symmetric monoidal category $\mathcal{E}$. In particular,
the strong monoidal functor $Set\rightarrow\mathcal{E}$, which sends
a set $S$ to the $S$-fold coproduct of copies of the unit $I$ of
$\mathcal{E}$, maps every operad $\mathcal{P}$ over $C$ in $Set$
to an operad in $\mathcal{E}$, which we denote by $\mathcal{P}_{\mathcal{E}}$,
or sometimes again by $\mathcal{P}$. 

If $\mathcal{P}$ is an operad over $C$ and $f:D\rightarrow C$ is
a map of sets, then there is an evident induced operad $f^{*}(\mathbb{\mathcal{P}})$
over $D$, given by \[
f^{*}(\mathcal{P})(d_{1},\cdots,d_{n};d)=\mathcal{P}(fd_{1},\cdots,fd_{n};fd).\]
 If $\mathcal{P}$ and $\mathcal{Q}$ are operads, a \emph{map} $\xymatrix{\mathbb{\mathcal{Q}}\ar[r]^{f} & \mathbb{\mathcal{P}}}
$ is given by a map of sets $f:D\rightarrow C$, and for each $d_{1},\cdots,d_{n},d$
a map \[
f_{d_{1},\cdots,d_{n},d}:\mathcal{Q}(d_{1},\cdots,d_{n};d)\rightarrow\mathcal{P}(f(d_{1}),\cdots,f(d_{n});f(d))\]
 which commutes with all the operations and the $\Sigma_{n}$-actions.
If $D=C$ and $f:D\rightarrow C$ is the identity, we will call $f$
a map of operads \emph{over} $C$. For a fixed symmetric monoidal
category $\mathcal{E}$, we denote by $Operad(\mathcal{E})$ the category
of all coloured operads in $\mathcal{E}$. When $\mathcal{E}=Set$
we will simply write $Operad$ instead of $Operad(Set)$. \\
~~~~~

\begin{example}
Let $\mathcal{E}$ be a symmetric monoidal category. Then $\mathcal{E}$
gives rise to a coloured operad $\underline{\mathcal{E}},$ whose
colours are the objects of $\mathcal{E}$. For a sequence $X_{1},\cdots,X_{n},X$
of such objects, $\underline{\mathcal{E}}(X_{1},\cdots,X_{n};X)$
is the set of arrows $X_{1}\otimes\cdots\otimes X_{n}\rightarrow X$
in $\mathcal{E}$. If $\mathcal{E}$ is a symmetric \emph{closed}
monoidal category, then $\mathcal{E}$ may be viewed as an operad
$\underline{\underline{\mathcal{E}}}$ in $\mathcal{E}$, with the
objects of $\mathcal{E}$ as colours again, and with $\underline{\underline{\mathcal{E}}}(X_{1},\cdots,X_{n};X)$
the internal Hom-object $\underline{Hom}_{\mathcal{E}}(X_{1}\otimes\cdots\otimes X_{n},X)$. 

Note that, in general, the objects of $\mathcal{E}$ form a proper
class and not a set. However, in this paper, we will largely ignore
such set-theoretic issues, and interpret ''small'' or ''set''
in terms of a suitable universe. In this context, let us point out
that for any \emph{set} $S$ of objects of $\mathcal{E}$, there are
operads $\underline{\mathcal{E}}_{S}$ and $\underline{\underline{\mathcal{E}}}_{S}$
obtained by restricting $\underline{\mathcal{E}}$ and $\underline{\underline{\mathcal{E}}}$
to the colours in S (If $i:S\rightarrow Objects(\mathcal{E})$ is
the inclusion, then $\underline{\mathcal{E}}_{S}=i^{*}(\underline{\mathcal{E}}),$
etc). In general, we will often identify a monoidal category with
the corresponding operad; and simply write $\mathcal{E}$ for $\underline{\mathcal{E}}$
or $\underline{\underline{\mathcal{E}}}$.
\end{example}

\begin{example}
Any category $\mathcal{C}$ can be considered as an operad $\mathcal{P}_{\mathcal{C}}$
in the following way. The colours of $\mathcal{P}_{\mathcal{C}}$
are the objects of $\mathcal{C}$, and for any sequence of colours
$c_{1},\cdots,c_{n},c$ we set \[
\mathcal{P}_{\mathcal{C}}(c_{1},\cdots,c_{n};c)=\left\{ \begin{array}{cc}
\mathcal{C}(c_{1},c), & \textrm{if }n=1\\
\phi, & \textrm{if }n\ne1\end{array}\right.\]
the compositions and units are as in $\mathcal{C}$ and the symmetric
actions are all trivial. In this way we obtain a functor $j_{!}:Cat\rightarrow Operad$
from the category $Cat$ of small categories to the category of operads.
This functor has an evident right adjoint $j^{*}:Operad\rightarrow Cat$,
sending an operad $\mathcal{P}$ to the category given by the colours
and unary operations of $\mathcal{P}$. In exactly the same way, any
$\mathcal{E}$-enriched category can be seen as an operad in $\mathcal{E}$
and we thus obtain adjoint functors $\xymatrix{Cat(\mathcal{E})\ar@<2pt>[r]^{j_{!}\,\,\,\,\,\,} & Operad(\mathcal{E}).\ar@<2pt>[l]^{j^{*}\,\,\,\,\,\,}}
$
\end{example}
\begin{rem}
There is also the notion of a non-symmetric (also called planar) operad.
A planar operad is exactly the same structure as an operad except
that there are no symmetric actions involved. The resulting category
of planar operads with their obvious notion of maps is denoted by
$Operad_{\pi}(\mathcal{E})$. There is an evident forgetful functor
$Operad(\mathcal{E})\rightarrow Operad(\mathcal{E})_{\pi}$ which
maps an operad to the same operad with the symmetric actions forgotten.
This functor has a left adjoint $Symm:Operad(\mathcal{E})_{\pi}\rightarrow Operad(\mathcal{E})$,
which we call the symmetrization functor. This functor is useful in
the construction of operads, since sometimes it is easier to directly
describe the non-symmetric operad whose algebras are the desired structures
in a given context. 
\end{rem}
\begin{example}
Let $S$ be a set. We describe now a planar operad $\mathcal{B}_{S}$
whose algebras are categories having $S$ as set of objects. The set
of colours of $\mathcal{B}_{S}$ is $S\times S$, and for any sequence
of colours of the form $(s_{1},s_{2}),(s_{2},s_{3}),\cdots,(s_{n-1},s_{n})$
there is exactly one operation in $\mathcal{B}_{S}((s_{1},s_{2}),\cdots,(s_{n-1},s_{n});(s_{1},s_{n}))$.
There are no other operations except those just given, which then
completely determine the operadic structure. We thus have a planar
operad in $Set$ whose symmetrization we denote by $\mathcal{A}_{S}$.
For any cocomplete monoidal category $\mathcal{E}$ we obtain an operad
in $\mathcal{E}$ (still) denoted $\mathcal{A}_{S}$ which is the
image of the original $\mathcal{A}_{S}$ under the functor $Operad(Set)\rightarrow Operad(\mathcal{E})$
described above. It is easy to verify that an $\mathcal{A}_{S}$-algebra
in $\mathcal{E}$ is the same as an $\mathcal{E}$-enriched category
having $S$ as set of objects. 
\end{example}
We refer the reader to \cite{BM col operads} for more examples of
coloured operads.

\section{A category of trees}

The trees we will consider are finite, non-empty (non-planar) trees with a designated root. As is common in the theory of operads
\cite{GetzKap,GinzKap,MSS}
we allow some edges to have a vertex only on one side. These edges are called \emph{outer} (or external) edges, while those having vertices
on both sides are called \emph{inner} (or internal) edges.
By a designated root we mean a choice of one of the outer edges. The
root defines an up-down direction in the tree (towards the root) and
thus each vertex has a number of incoming edges (the number is the
\emph{valence} of the vertex) and one edge going out of it. We also
allow vertices of valence 0. For example, the tree 
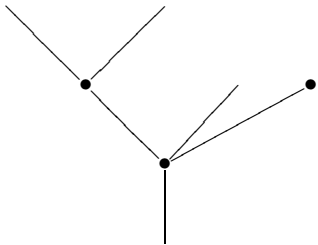
\[
\xymatrix{*{\,}\ar@{-}[dr] &  & *{\,}\ar@{-}[dl]\\
 & *{\bullet}\ar@{-}[dr] &  & *{\,}\ar@{-}[dl] & *{\bullet}\ar@{-}[dll]\\
 &  & *{\bullet}\ar@{-}[d]\\
 &  & *{\,}}
\]
has three vertices, of valence 2, 3, and 0, and three input edges.
A tree with no vertices \[
\xymatrix{*{\,}\ar@{-}[d]\\
*{\,}}
\]
whose input edge ($e$ say) coincides with its output edge will be
denoted by $\eta_{e}$.

When we draw a tree we will always put the root at the bottom. One
drawback of drawing a tree on the plane is that it immediately becomes
a planar tree; we thus have many different 'pictures' for the same
tree. For instance the two trees
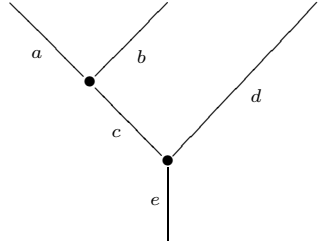
\[
\xymatrix{*{}\ar@{-}[dr]_{a} &  & *{}\ar@{-}[dl]^{b} &  & *{}\ar@{-}[ddll]^{d}\\
 & *{\bullet}\ar@{-}[dr]_{c}\\
 &  & *{\bullet}\ar@{-}[d]_{e}\\
 &  & *{}}
\]
 and
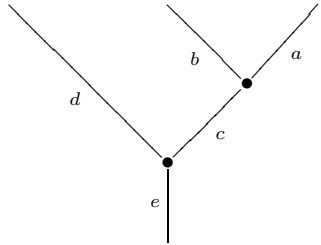
\[
\xymatrix{*{}\ar@{-}[ddrr]_{d} &  & *{}\ar@{-}[dr]_{b} &  & *{}\ar@{-}[dl]^{a}\\
 &  &  & *{\bullet}\ar@{-}[dl]^{c}\\
 &  & *{\bullet}\ar@{-}[d]_{e}\\
 &  & *{}}
\]
 are different planar representations of the same tree. 

Any tree $T$ can be viewed as generating an operad $\Omega(T)$,
whose colours are the edges of the tree, while the vertices of the
tree are the generators of the operations. More explicitly, if we
choose a planar representation of $T$ then each vertex $v$ with
input edges $e_{1},\cdots,e_{n}$ and output edge $e$ defines an
operation $v\in\Omega(T)(e_{1},\cdots,e_{n};e)$. The other operations
are the unit operations and the operations obtained by compositions
and by permutations, so as to obtain an operad in which every Hom
set has at most one object. For example, in the same tree $T$ pictured
above, let us name the edges and vertices $a,b,\cdots,f$ and $r,v,w.$ 

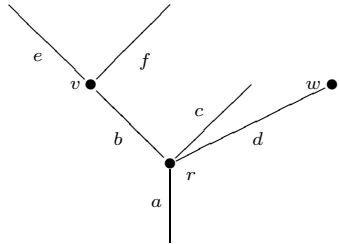
\[
\xymatrix{*{\,}\ar@{-}[dr]_{e} &  & *{\,}\ar@{-}[dl]^{f}\\
\,\ar@{}[r]|{\,\,\,\,\,\,\,\,\,\,\,\,\,\, v} & *{\bullet}\ar@{-}[dr]_{b} &  & *{\,}\ar@{-}[dl]_{c}\ar@{}[r]|{\,\,\,\,\,\,\,\,\,\,\,\, w} & *{\bullet}\ar@{-}[dll]^{d}\\
 &  & *{\bullet}\ar@{-}[d]_{a} & \,\ar@{}[l]^{r\,\,\,\,\,\,\,\,\,\,\,}\\
 &  & *{\,}}
\]
 Then $v\in\Omega(T)(e,f;b),\, w\in\Omega(T)(\,;d)$ and $r\in\Omega(b,c,d;a))$
are the generators, while the other operations are the units $1_{a},1_{b},1_{c}\cdots 1_{f} $,
the operations obtained by compositions $r\circ_{1}v\in\Omega(T)(e,f,c,d;a)$,
$r\circ_{3}w\in\Omega(T)(b,c;a)$ and $r(v,1_{c},w)=(r\circ_{1}v)\circ_{4}w=(r\circ_{3}w)\circ_{1}v\in\Omega(T)(e,f,c;a)$,
and permutations of these. This is a complete description of the operad
$\Omega(T)$. 

Viewing trees as coloured operads as above enables us to define the
category $\Omega$, whose objects are trees, and whose arrows $T\rightarrow T'$
are operad maps $\Omega(T)\rightarrow\Omega(T')$. The category $\Omega$
extends the simplicial category $\Delta$. Indeed, any $n\ge0$ defines
a linear tree
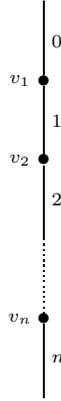
\[
\xymatrix{ & *{\,}\ar@{-}[d]^{0}\\
\,\ar@{}[r]|{\,\,\,\,\,\,\,\,\, v_{1}} & *{\bullet}\ar@{-}[d]^{1}\\
\,\ar@{}[r]|{\,\,\,\,\,\,\,\,\, v_{2}} & *{\bullet}\ar@{-}[d]^{2}\\
 & *{\,}\ar@{..}[d]\\
\,\ar@{}[r]|{\,\,\,\,\,\,\,\,\, v_{n}} & *{\bullet}\ar@{-}[d]^{n}\\
 & *{\,}}
\]
 on $n+1$ edges and $n$ vertices $v_{1},\cdots,v_{n}$. We denote
this tree by $[n]$. Any order preserving map $\{0,\cdots,n\}\rightarrow\{0,\cdots,m\}$
defines an arrow $[n]\rightarrow[m]$ in the category $\Omega$. In
this way, we obtain an embedding\[
\xymatrix{\Delta\ar[r]^{i} & \Omega}
\]
 This embedding is fully faithful. Moreover, it describes $\Delta$
as a sieve (or ideal) in $\Omega$, in the sense that for any arrow
$S\rightarrow T$ in $\Omega$, if $T$ is linear then so is $S$. 

With a tree $T$ one can associate certain maps in $\Omega$ as follows.
If $b$ is an inner edge in $T$, let $T/b$ be the tree obtained
from $T$ by contracting $b$. Then there is a natural map $\partial_{b}:T/b\rightarrow T$
in $\Omega$, called the \emph{inner face} map associated with
$b$, which locally in the tree looks like this:\\

$
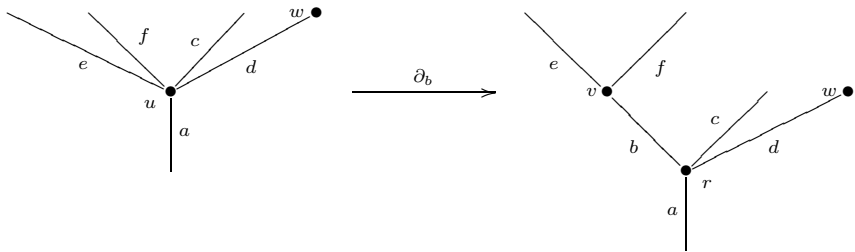
\xymatrix{*{\,}\ar@{-}[rrd]_{e} & *{\,}\ar@{-}[rd]^{f} &  & *{\,}\ar@{-}[dl]_{c}\ar@{}[r]|{\,\,\,\,\,\,\,\,\, w} & *{\bullet}\ar@{-}[lld]^{d}\\
 & \,\ar@{}[r]_{\,\,\,\,\,\,\,\,\,\,\, u} & *{\bullet}\ar@{-}[d]^{a}\\
 &  & *{\,}}
$ $\xymatrix{\\\,\ar[rr]^{\partial_{b}} &  & \,}
$$\xymatrix{*{\,}\ar@{-}[dr]_{e} &  & *{\,}\ar@{-}[dl]^{f}\\
\,\ar@{}[r]|{\,\,\,\,\,\,\,\,\,\,\,\,\,\, v} & *{\bullet}\ar@{-}[dr]_{b} &  & *{\,}\ar@{-}[dl]_{c}\ar@{}[r]|{\,\,\,\,\,\,\,\,\,\,\,\, w} & *{\bullet}\ar@{-}[dll]^{d}\\
 &  & *{\bullet}\ar@{-}[d]_{a} & \,\ar@{}[l]^{r\,\,\,\,\,\,\,\,\,\,\,}\\
 &  & *{\,}}
$ 

Let $v$ be a vertex in $T$ with the property that all but one of
the edges incident to $v$ are outer. We call such a vertex an
outer cluster. Let $T/v$ be the tree obtained from $T$ by removing
the vertex $v$ and all of the outer edges incident to it. Then
there is a map $\partial_{v}:T/v\rightarrow T$ in $\Omega$ called
the \emph{outer face} associated with $v$. For example, the maps:

~~~~~~~~~~~~
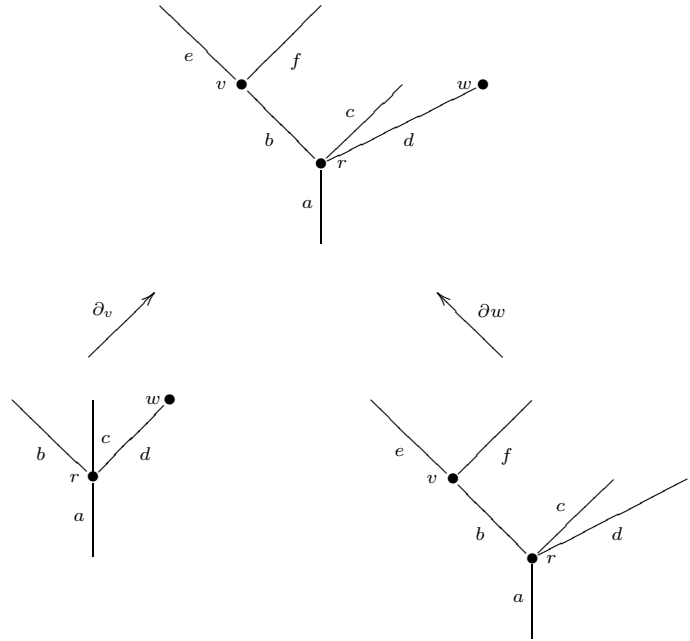
\[
\xymatrix{*{\,}\ar@{-}[dr]_{e} &  & *{\,}\ar@{-}[dl]^{f}\\
\,\ar@{}[r]|{\,\,\,\,\,\,\,\,\,\,\, v} & *{\bullet}\ar@{-}[dr]_{b} &  & *{\,}\ar@{-}[dl]_{c}\ar@{}[r]|{\,\,\,\,\,\,\,\,\,\,\,\, w} & *{\bullet}\ar@{-}[dll]^{d}\\
 &  & *{\bullet}\ar@{-}[d]_{a} & \,\ar@{}[l]|{r\,\,\,\,\,\,\,\,\,\,\,}\\
 &  & *{\,}}
\]
$ $$ $ $ $\\
~$\textrm{ $ $ }$~~~~~~~~~~~~~~~~~~~~~~~~$\xymatrix{ & \,\\
\ar[ur]^{\partial_{v}}}
$~~~~~~~~~~~~~~~~~~~~~~~~~~~~$\xymatrix{\,\\
 & \ar[ul]_{\partial w}}
$\\
\\
~~$\,\,\textrm{ $ $ }$~~~~~~~~~~~~~~~~$\xymatrix{*{\,}\ar@{-}[dr]_{b} & *{\,}\ar@{-}[d]^{c}\ar@{}[r]|{\,\,\,\,\,\,\,\,\,\,\,\, w} & *{\bullet}\ar@{-}[dl]^{d}\\
\,\ar@{}[r]|{\,\,\,\,\,\,\,\,\,\,\,\, r} & *{\bullet}\ar@{-}[d]_{a}\\
 & *{\,}}
$ ~~~~~~~~~~~~~~~~~~~~$\xymatrix{*{\,}\ar@{-}[dr]_{e} &  & *{\,}\ar@{-}[dl]^{f}\\
\,\ar@{}[r]|{\,\,\,\,\,\,\,\,\,\,\, v} & *{\bullet}\ar@{-}[dr]_{b} &  & *{\,}\ar@{-}[dl]_{c} & *{}\ar@{-}[dll]^{d}\\
 &  & *{\bullet}\ar@{-}[d]_{a} & \,\ar@{}[l]|{r\,\,\,\,\,\,\,\,\,\,\,\,}\\
 &  & *{\,}}
$ \\
are two outer faces. We will use the term \emph{face} map to refer
to an inner or outer face map. One more type of map is a map that can be associated with 
a unary vertex $v$ in $T$ as follows. Let $T/v$ be the tree obtained from $T$ by removing the vertex $v$ and merging the two edges
incident to it into one edge $e$. 
Then there is a map $\sigma_{v}:T\rightarrow T/v$
in $\Omega$ called the \emph{degeneracy map} associated with $v$,
which sends the vertex $v$ to the identity $1_{e},$ and
which can be pictured like this:\\

$
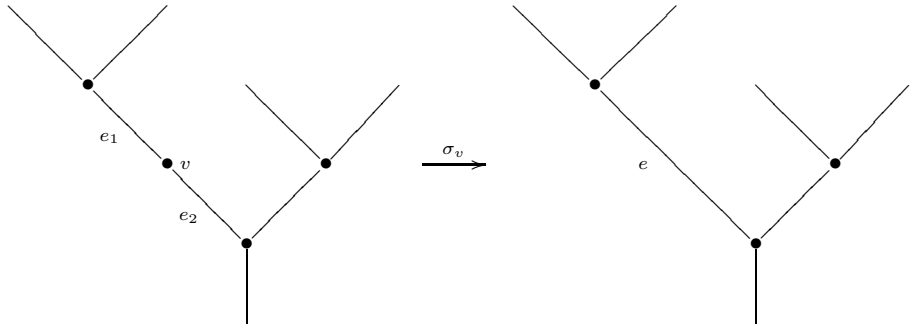
\xymatrix{*{\,}\ar@{-}[dr] &  & *{\,}\ar@{-}[dl]\\
 & *{\bullet}\ar@{-}[dr]_{e_{1}} &  & *{\,}\ar@{-}[dr] &  & *{\,}\ar@{-}[dl]\\
 &  & *{\bullet}\ar@{-}[dr]_{e_{2}}\ar@{}|{\,\,\,\,\,\,\,\,\,\, v} &  & *{\bullet}\ar@{-}[dl]\\
 &  &  & *{\bullet}\ar@{-}[d]\\
 &  &  & *{\,}}
$$\xymatrix{\\\\\,\ar[r]^{\sigma_{v}} & \,}
$$\xymatrix{*{\,}\ar@{-}[dr] &  & *{\,}\ar@{-}[dl]\\
 & *{\bullet}\ar@{-}[ddrr] &  & *{\,}\ar@{-}[dr] &  & *{\,}\ar@{-}[dl]\\
 & \,\ar@{}[r]|{\,\,\,\, e} &  &  & *{\bullet}\ar@{-}[dl]\\
 &  &  & *{\bullet}\ar@{-}[d]\\
 &  &  & *{\,}}
$ 

The following lemma is the generalization to $\Omega$ of the well
known fact that in $\Delta$ each arrow can be written as a composition
of degeneracy maps followed by face maps. We omit the proof.

\begin{lem}
Any arrow $f:A\rightarrow B$ in $\Omega$ decomposes as \[
\xymatrix{A\ar[r]^{f}\ar[d]^{\sigma} & B\\
A'\ar[r]^{\varphi} & B'\ar_{\delta}[u]}
\]
 where $\sigma:A\rightarrow A'$ is a composition of degeneracy maps,
$\varphi:A'\rightarrow B'$ is an isomorphism, and $\delta:B'\rightarrow B$
is a composition of face maps. 
\end{lem}

\section{Dendroidal sets}

We now define the category $dSet$ of dendroidal sets and discuss its relation to the category $sSet$ of simplicial sets. 

\begin{defn}
A \emph{dendroidal set} is a functor $\Omega^{op}\rightarrow Set$.
A map between dendroidal sets is a natural transformation. The category
of \emph{dendroidal sets} thus defined is denoted \emph{dSet}.
\end{defn}
Thus, a dendroidal set $X$ is given by a set $X_{T}$ for each tree
$T$, and a map $\alpha^{*}:X_{T}\rightarrow X_{S}$ for each map
of trees (arrow in $\Omega$) $\alpha:S\rightarrow T$; and these
maps have to be functorial in $\alpha$, in the sense that $id^{*}=id$
and $(\alpha\beta)^{*}=\beta^{*}\alpha^{*}$ for $\xymatrix{R\ar[r]^{\beta} & S\ar[r]^{\alpha} & T}
$ in $\Omega$. A morphism $\xymatrix{Y\ar[r]^{f} & X}
$ of dendroidal sets is given by maps (all denoted) $f:Y_{T}\rightarrow X_{T}$
for each tree $T$, commuting with the structure maps (i.e., $f(\alpha^{*}y)=\alpha^{*}f(y)$
for any $y\in Y_{T}$ and any $\alpha:S\rightarrow T$). An element
of $X_{T}$ is called a \emph{dendrex} (plural \emph{dendrices}) of
\emph{shape} $T$ (This terminology is analogous to simplex, simplices).
The dendrices of shape $\eta$ will be referred to as \emph{vertices}. As for 
simplicial sets, we call a dendrex $x\in X_{T}$ \emph{degenerate} if there exists a 
degeneracy $\sigma :T\twoheadrightarrow S$ and a dendrex $y\in X_{S}$ with $\sigma ^{*}(y)=x$.

Every tree $T$ defines a representable dendroidal set $\Omega[T]$
as follows:\[
\Omega[T]_{S}=\Omega(S,T).\]
By the Yoneda Lemma each dendrex $x$ of shape $T$ in a dendroidal
set $X$ corresponds bijectively to a map $\hat{x}:\Omega[T]\rightarrow X$
of dendroidal sets. If $\partial_{x}:T\rightarrow R$ is a face map
associated to an inner edge or an outer cluster $x$ we use
the same notation $\partial_{x}:\Omega[T]\rightarrow\Omega[R]$ for
the induced map of dendroidal sets. 

The inclusion functor $i:\Delta\rightarrow\Omega$ defines an obvious
restriction functor\[
i^{*}:dSet\rightarrow sSet.\]
 This functor has both a left adjoint $i_{!}$ and a right adjoint
$i_{*}$, given by left and right Kan extension. The functor $i_{!}:sSet\rightarrow dSet$
is ''extension by zero'', \[
i_{!}(X)_{T}=\left\{ \begin{array}{cc}
X_{n}, & \textrm{if T is linear with $n$ vertices}\\
\phi, & \textrm{otherwise}\end{array}\right.\]
 (This is clear from the fact that $\Delta\subseteq\Omega$ is a sieve).
It follows that $i_{!}$ is full and faithful, and that $i^{*}i_{!}$
is the identity functor on simplicial sets. The pair $(i^{*},i_{*})$
defines a morphism of toposes $i:sSet\rightarrow dSet$, which is
in fact an \emph{open} \emph{embedding}.

\begin{example}
If $\mathcal{P}$ is an operad, then the \emph{dendroidal nerve} of
$\mathcal{P}$ is the dendroidal set $N_{d}(\mathcal{P})$ given by\[
N_{d}(\mathcal{P})_{T}=Hom_{Operad}(\Omega(T),\mathcal{P}),\]
 
\end{example}
This construction defines a fully faithful functor \[
N_{d}:Operad\rightarrow dSet,\]
 which has various nice properties as we will see. As already noted,
any monoidal category $\mathcal{E}$ defines an operad $\underline{\mathcal{E}}$.
The corresponding dendroidal set $N_{d}(\underline{\mathcal{E}})$
will simply be written $N_{d}(\mathcal{E})$ and will be called the
\emph{dendroidal nerve} of $\mathcal{E}$. Note that this extends
the usual (simplicial) nerve of $\mathcal{E}$, in the sense that
\[
i^{*}(N_{d}\mathcal{E})=N(\mathcal{E}).\]
The functor $N_{d}:Operad\rightarrow dSet$ has a left adjoint \[
\tau_{d}:dSet\rightarrow Operad\]
defined by Kan extension. For a dendroidal set $X$, we refer to $\tau_{d}(X)$
as the \emph{operad generated by $X$}. This functor $\tau_{d}$ extends
the functor $\tau$ from simplicial sets to categories, left adjoint
to $N:Cat\rightarrow sSet$. In particular, we obtain a diagram of
functors \[
\xymatrix{sSet\ar@<2pt>[r]^{i_{!}}\ar@<-2pt>[d]_{\tau} & dSet\ar@<2pt>[l]^{i^{*}}\ar@<-2pt>[d]_{\tau_{d}}\\
Cat\ar@<-2pt>[u]_{N}\ar@<2pt>[r]^{j_{!}} & Operad\ar@<-2pt>[u]_{N_{d}}\ar@<2pt>[l]^{j^{*}}}
\]
(with left adjoints on the top or on the left), in which the following
commutation relations hold up to natural isomorphisms\[
\tau N=id,\,\,\,\tau_{d}N_{d}=id,\,\,\, i^{*}i_{!}=id,\,\,\, j^{*}j_{!}=id\]
 and \[
j_{!}\tau=\tau_{d}i_{!},\,\,\, Nj^{*}=i^{*}N_{d},\,\,\, i_{!}N=N_{d}j_{!}.\]
(The canonical map $\tau i^{*}(X)\rightarrow j^{*}\tau _{d}(X)$ is in general not an isomorphism.)

\begin{rem}
For an arbitrary category $\mathcal{E}$, one can also consider the
category $d\mathcal{E}$ of dendroidal objects in $\mathcal{E}$,
i.e., contra-variant functors from $\Omega$ to $\mathcal{E}$. In
particular, if one takes for $\mathcal{E}$ the category $Top$ of
compactly generated topological spaces, one obtains in this way the
category $dTop$ of dendroidal spaces. Many constructions extend to
this more general context. For example, if $\mathcal{P}$ is a topological
operad, its dendroidal nerve $N_{d}(\mathcal{P})$ is naturally a
dendroidal space, with the special property that its space $N_{d}(\mathcal{P})_{\eta}$
of vertices is discrete. Conversely, from such a dendroidal space
$X$ with this property, one can construct a topological operad, $\tau_{d}(X)$. 
\end{rem}

\subsection{Diagrams of dendroidal sets}

If $X:\mathbb{S}^{op}\rightarrow sSet$ is a diagram of simplicial
sets (contravariantly) indexed by a small category $\mathbb{S}$,
one can construct a ''total'' simplicial set $\int_{\mathbb{S}}X$ as follows.
An $n$-simplex of $\int_{\mathbb{S}}X$ is a pair $(s,x)$ where $s=(\xymatrix{s_{0}\ar[r]^{\alpha_{1}} & s_{1}\ar[r] & \cdots\ar[r]^{\alpha_{n}} & s_{n}}
)$ is an $n$-simplex in the nerve of $\mathbb{S}$, and $x$ is a function
assigning to each map $u:[k]\rightarrow[n]$ in $\Delta$ a $k$-simplex
$x_{k}$ in $X(s_{u(0)})$, functorial in the following way. If $w=uv$
:$\xymatrix{[l]\ar[r]^{v} & [k]\ar[r]^{u} & [n]}
$, then $u(0)\le w(0)$ so there is a composition of $\alpha_{i}$'s
from $s_{u(0)}$ to $s_{w(0)}$ in $\mathbb{S}$, denoted $\alpha_{w,u}=\alpha_{w(0)}\circ\alpha_{w(0)-1}\circ\cdots\circ\alpha_{u(0)+1}$.
Then the functorial condition on the $x_{\alpha}$'s is \[
\alpha_{w,u}^{*}(x_{w})=v^{*}(x_{u})\]
 (here $\alpha_{w,u}^{*}:X(s_{w(0)})\rightarrow X(s_{u(0)})$, and
this is an identity between $l$-simplices in $X(s_{u(0)})$ ). Notice
that in the special case where we start with a diagram $\mathbb{C}:\mathbb{S}^{op}\rightarrow Cat$
of small categories, the diagram $Nerve(\mathbb{C}):\mathbb{S}^{op}\rightarrow sSet$
satisfies the identity\[
\int_{\mathbb{S}}(Nerve(\mathbb{C}))=Nerve(\int_{\mathbb{S}}\mathbb{C})\]
 where $\int_{\mathbb{S}}\mathbb{C}$ on the right is the Grothendieck
construction. 

We shall now give a similar construction for diagrams of dendroidal
sets. For this, we assume that the indexing category $\mathbb{S}$
has finite products. So, let $X:\mathbb{S}^{op}\rightarrow dSet$
be a diagram of dendroidal sets. We define a dendroidal set $\int_{\mathbb{S}}X$
as follows. For a tree $T$, an element of $\int_{\mathbb{S}}X_{T}$ is again
a pair $(t,x)$. Here $t\in N_{d}(\mathbb{S})_{T}$ is an element
of the dendroidal nerve of $\mathbb{S}$ (where $\mathbb{S}$ is viewed
as an operad via the cartesian structure). Such an element determines
an object $in(t)\in\mathbb{S}$, defined by $in(t)=t(e_{1})\times\cdots\times t(e_{n})$
where $e_{1},\cdots,e_{n}$ are the input edges of $T$ (in some fixed
arbitrary order). Note that for any arrow $u:S\rightarrow T$ in $\Omega$,
the dendrex $t$ determines a map $in(t)\rightarrow in(tu)$ in $\mathbb{S}$
(defined by projections, the maps given by $t$, and the coherence
maps in $\mathbb{S}$). Now $x$ is a function which assigns to each
such $u$ an element $x_{u}\in X(in(tu))_{S}$, functorial in the
following way: if $w=u\circ v$ as in $\xymatrix{R\ar[r]^{v} & S\ar[r]^{u} & T}
$ then there is an induced map $\xymatrix{in(tu)\ar[r]^{\alpha_{u,v}} & in(tw)}
$ in $\mathbb{S}$, and we require \[
\alpha_{u,v}^{*}(x_{v})=v^{*}(x_{u})\]

The set $\int_{\mathbb{S}}X_{T}$ of such pairs $(t,x)$ is contravariant in $T$,
and defines the dendroidal set $\int_{\mathbb{S}}X$.

Note that this construction for dendroidal sets truly extends the
one for simplicial sets, in the sense that for a diagram $X:\mathbb{S}^{op}\rightarrow sSet$
of simplicial sets where $\mathbb{S}$ is cartesian, there is a canonical
isomorphism\[
i_{!}\int_{\mathbb{S}}X=\int_{\mathbb{S}}i_{!}X.\]

\section{The tensor product of dendroidal sets}

Like any other category of presheaves of sets, the category $dSet$
has a closed cartesian structure. There is, however, another more
interesting monoidal structure on the category of dendroidal sets,
which we aim to describe in this section. To begin with, we will recall
the tensor product for operads from \cite{BV hom. inv. str.}.

\subsection{The Boardman-Vogt tensor product}

Let $\mathcal{P}$ be an operad in $Set$ over $C$ and $\mathcal{Q}$
one over $D$. Their tensor product $\mathcal{P}\otimes_{BV}\mathcal{Q}$
is an operad coloured by the product set $C\times D$. The operations
in $\mathcal{P}\otimes_{BV}\mathcal{Q}$ are generated by the following.
Any $p\in\mathcal{P}(c_{1},\cdots,c_{n};c)$ and any $d\in D$ define
an operation\[
p\otimes d\in\mathcal{P}\otimes_{BV}Q((c_{1},d),\cdots,(c_{n},d);(c,d)).\]
 These operations compose in $\mathcal{P}\otimes_{BV}\mathcal{Q}$
in a way to make $p\mapsto p\otimes d$ a map of operads. Similarly,
each operation $q\in\mathcal{Q}(d_{1},\cdots,d_{n},d)$ and each $c\in C$
define an operation \[
c\otimes q\in\mathcal{P}\otimes_{BV}\mathcal{Q}((c,d_{1}),\cdots(c,d_{n});(c,d)),\]
and these compose as in $Q$. Furthermore, the operations from $\mathcal{P}$
and $\mathcal{Q}$ distribute over each other, in the sense that for
$p\in\mathcal{P}(c_{1},\cdots,c_{n};c)$ and $q\in\mathcal{Q}(d_{1},\cdots,d_{m};d),$\[
\sigma_{n,m}^{*}((p\otimes d)(c_{1}\otimes q,\cdots,c_{n}\otimes q))=(c\otimes q)(p\otimes d_{1},\cdots,p\otimes d_{m})\]
 where $\sigma_{n,m}\in\Sigma_{n\cdot m}$ is the permutation described
as follows. Consider $\Sigma_{n\cdot m}$ as the set of bijections
of the set $\{0,1,\cdots,n\cdot m-1\}$. Each number in this set can
be written uniquely in the form $k\cdot n+j$ where $0\le k<m$ and
$0\le j<n$ as well as in the form $k\cdot m+j$ where $0\le k<n$
and $0\le j<m$. The permutation $\sigma_{n,m}$ is then defined by
$\sigma_{n,m}(k\cdot n+j)=j\cdot m+k$. This tensor product makes
the category of operads into a symmetric monoidal category.

This Boardman-Vogt tensor product preserves colimits in each variable
separately. In fact, there is a corresponding internal Hom, making
the category $Operad$ into a symmetric \emph{closed} monoidal category.
For two operads $\mathcal{P}$ and $\mathcal{Q}$ as above, $\underline{Hom}(\mathcal{P},\mathcal{Q})$
is the operad whose colours are the maps $\mathcal{P}\rightarrow\mathcal{Q}$,
and whose operations are suitably defined multi-natural transformations.
(Explicitly, for $\alpha_{1},\cdots,\alpha_{n},\beta:\mathcal{P}\rightarrow\mathcal{Q}$,
elements of $\underline{Hom}(\mathcal{P},\mathcal{Q})(\alpha_{1},\cdots,\alpha_{n};\beta)$
are maps $f$ assigning to each colour $c\in C$ of $\mathcal{P}$
an element $f_{c}\in\mathcal{Q}(\alpha_{1}c,\cdots,\alpha_{n}c;\beta c).$
These $f_{c}$ should be natural with respect to all operations in
$\mathcal{P}$. For example, if $p\in\mathcal{P}(c_{1},c_{2};c)$
is a binary operation, then $\beta(p)(f_{c_{1}},f_{c_{2}})\in\mathcal{Q}(\alpha_{1}c_{1},\cdots,\alpha_{n}c_{1},\alpha_{1}c_{2},\cdots,\alpha_{n}c_{2};\beta c)$
is the image under a suitable permutation of $f_{c}(\alpha_{1}(p),\cdots,\alpha_{n}(p))\in\mathcal{Q}(\alpha_{1}c_{1},\alpha_{1}c_{2},\cdots,\alpha_{n}c_{1},\alpha_{n}c_{2};\beta c)$).

For a symmetric monoidal category $\mathcal{E}$, the Boardman-Vogt
tensor product of coloured operads in $\mathcal{E}$ still makes sense
for \emph{Hopf} operads $\mathcal{P}$ and $\mathcal{Q}$.
For such operads, the categories $Alg_{\mathcal{E}}(\mathcal{P})$
and $Alg_{\mathcal{E}}(\mathcal{Q})$ are again symmetric monoidal,
and a $(\mathcal{P}\otimes_{BV}\mathcal{Q})$-algebra in $\mathcal{E}$
is the same thing as a $\mathcal{P}$-algebra in $Alg_{\mathcal{E}}(\mathcal{Q})$,
and is also the same thing as a $\mathcal{Q}$-algebra in $Alg_{\mathcal{E}}(\mathcal{P})$.

\subsection{The tensor product of dendroidal sets}

We now define a tensor product \[
\otimes:dSet\times dSet\rightarrow dSet\]
 which is to preserve colimits in each variable separately. Since
each dendroidal set is a colimit of representables, this tensor is
completely determined by its effect on representable dendroidal sets
$\Omega[S]$ and $\Omega[T]$, which we define as \[
\Omega[S]\otimes\Omega[T]=N_{d}(\Omega(S)\otimes_{BV}\Omega(T)),\]
 i.e., as the dendroidal nerve of the Boardman-Vogt tensor product
of the operads $\Omega(S)$ and $\Omega(T)$. It follows by general
category theory \cite{Day,Kelly ECT} that there exists an internal
Hom for this tensor, defined for two dendroidal sets $X$ and $Y$
and an object $T$ of $\Omega$ by\[
\underline{Hom}(X,Y)_{T}=Hom_{dSet}(\Omega[T]\otimes X,Y)\]

We summarise this discussion in the following proposition:

\begin{prop}
There exists a unique (up to natural isomorphism) symmetric closed
monoidal structure on $dSet$, with the property that there is a natural
isomorphism $\Omega[S]\otimes\Omega[T]\cong N_{d}(\Omega(S)\otimes_{BV}\Omega(T))$
for any two objects $S,T$ of $\Omega$.
\end{prop}
More generally, for suitable symmetric monoidal categories $\mathcal{E}$,
there is such a monoidal structure on the category $d\mathcal{E}$
of dendroidal objects. See the appendix for a discussion of dendroidal
objects.

We mention some basic properties of the tensor product on $dSet$,
in relation to the tensor product of operads, and to the product of
simplicial sets. 

\begin{prop}
The following properties hold.

(i) For any two dendroidal sets $X$ and $Y$, there is a natural
isomorphism \[
\tau_{d}(X\otimes Y)\cong\tau_{d}(X)\otimes_{BV}\tau_{d}(Y).\]

(ii) For any two operads $\mathcal{P}$ and $\mathcal{Q}$, there
is a natural isomorphism\[
\tau_{d}(N_{d}(\mathcal{P})\otimes N_{d}(\mathcal{Q}))\cong\mathcal{P}\otimes_{BV}\mathcal{Q}.\]

\end{prop}
\begin{proof}
It suffices to check $(i)$ for representable $X$ and $Y$, in which
case it follows from the identity $\tau_{d}N_{d}\cong id.$ By the
same identity, $(ii)$ follows from $(i)$.
\end{proof}
\begin{prop}
For any two simplicial sets $X$ and $Y$, and any dendroidal set
$D$, there are natural isomorphisms\\

(i) $i_{!}(X\times Y)\cong i_{!}(X)\otimes i_{!}(Y)$,\\
~

(ii) $i^{*}\underline{Hom}(i_{!}(X),D)\cong i^{*}(D)^{X}$,

~

(iii) $i^{*}\underline{Hom}(i_{!}(X),i_{!}(Y))\cong Y^{X}$.
\end{prop}
\begin{proof}
The isomorphisms of type $(ii)$ and $(iii)$ are deduced from those
of type $(i)$, using the fact that $i_{!}$ is fully faithful. For
$(i)$, it suffices again to check this for representable simplicial
sets $\Delta[n]$ and $\Delta[m]$. Observe first that, more generally,
for any two small categories $\mathbb{C}$ and $\mathbb{B}$, \[
\quad\quad\quad\quad\quad\quad\quad\quad\quad\quad j_{!}(\mathbb{C}\times\mathbb{B})\cong j_{!}(\mathbb{C})\otimes_{BV}j_{!}(\mathbb{B})\quad\quad\quad\quad\quad\quad\quad\quad\quad\quad\,\,\,\,\,\,\,\,\,\,(1)\]
 This holds in particular for the linear orders $[n]$ and $[m]$
viewed as categories, so \begin{eqnarray*}
i_{!}(\Delta[n]\times\Delta[m]) & \cong & i_{!}(N([n])\times N([m]))\\
 & \cong & i_{!}(N([n]\times[m]))\\
 & \cong & N_{d}j_{!}([n]\times[m])\\
 & \cong & N_{d}(j_{!}[n]\otimes_{BV}j_{!}[m])\,\,\,\,\textrm{(by (1))}\\
 & \cong & N_{d}(\Omega(n)\otimes_{BV}\Omega(m))\\
 & \cong & \Omega[n]\otimes\Omega[m]\\
 & \cong & i_{!}(\Delta[n])\otimes i_{!}(\Delta[m]).\end{eqnarray*}
This shows that $(i)$ holds for representables $\Delta[n]$ and $\Delta[m]$;
as said, this completes the proof. 
\end{proof}

\section{The homotopy coherent nerve}

We begin by recalling the Boardman-Vogt resolution of operads \cite{BV hom. inv. str.}
and its generalization \cite{BM W-cons}. 

Let $\mathcal{P}=(C,P)$ be an operad in the category of compactly
generated topological spaces, and let $H=[0,1]$ be the unit interval.
One can construct a (cofibrant) resolution $W(\mathcal{P})\rightarrow\mathcal{P}$
as follows. $W(\mathcal{P})$ is again an operad coloured by $C$.
The space $W(\mathcal{P})(c_{1},\cdots,c_{n};c)$ is a quotient of
a space of labelled planar trees. The edges of such a tree are labeled
by elements of $C$, where in particular the input edges are labelled
by the given $c_{1},\cdots,c_{n}$ and the output by $c$. Moreover,
the inner edges carry a label $t\in H$ (a ''length''), and each
vertex $v$ with input edges labelled $b_{1},\cdots,b_{n}\in C$ (in
the planar order) and output edge labelled $b\in C$, is labelled
by an element $p\in\mathcal{P}(b_{1},\cdots,b_{n};b)$. For example,
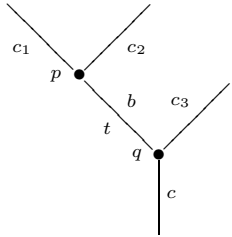
\[
\xymatrix{\ar@{-}[dr]_{c_{1}} &  & \ar@{-}[dl]^{c_{2}}\\
\ar@{}[r]|{\,\,\,\,\,\,\,\,\, p} & *{\bullet}\ar@{-}[dr]^{b}\ar@{}[dr]_{t} &  & \ar@{-}[dl]_{c_{3}}\\
 & \ar@{}[r]|{\,\,\,\,\,\,\,\,\,\, q} & *{\bullet}\ar@{-}[d]^{c}\\
 &  & *{\,}}
\]
 where $p\in\mathcal{P}(c_{1},c_{2};b)$, $q\in\mathcal{P}(b,c_{3};c)$,
$t\in[0,1]$. There is a natural (product) topology on these trees,
coming from the topology on $\mathcal{P}$ and that on $H$. The space
$W(\mathbb{\mathcal{P}})(c_{1},\cdots,c_{n};c)$ is now the quotient
space, obtained (by identifying isomorphic planar trees with the same
labeling and) the following two relations:

(i) Vertices labeled by an identity can be deleted, taking the maximum
of the two adjacent lengths (or forgetting the lengths altogether
if one of the adjacent edges is outer) 
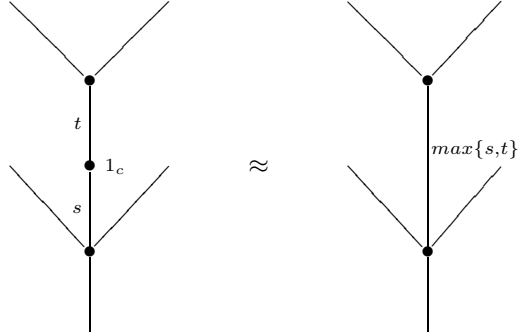
\[
\xymatrix{*{\,}\ar@{-}[dr] &  & *{\,}\ar@{-}[dl] &  & *{\,}\ar@{-}[dr] &  & *{\,}\ar@{-}[dl]\\
 & *{\bullet}\ar@{-}[d]_{t} &  &  &  & *{\bullet}\ar@{-}[dd]\\
*{\,}\ar@{-}[dr] & *{\bullet}\ar@{-}[d]_{s} & *{\,}\ar@{-}[dl]\ar@{}[l]|{1_{c}\,\,\,\,\,\,\,\,} & \approx & *{\,}\ar@{-}[dr] &  & *{\,}\ar@{-}[dl]\ar@{}[l]_{\,\,\,\,\,\, max\{ s,t\}}\\
 & *{\bullet}\ar@{-}[d] &  &  &  & *{\bullet}\ar@{-}[d]\\
 & *{\,} &  &  &  & *{\,}}
\]
$ $ 

(ii) Edges of length zero can be contracted, using the operad composition
of $\mathcal{P}$: 
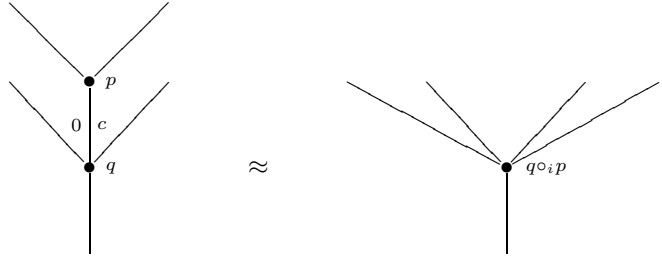
\[
\xymatrix{*{\,}\ar@{-}[dr] &  & *{\,}\ar@{-}[dl]\\
*{\,}\ar@{-}[dr] & *{\bullet}\ar@{-}[d]_{0}^{c} & *{\,}\ar@{-}[dl]\ar@{}[l]|{p\,\,\,\,\,\,\,\,\,\,} &  & *{\,}\ar@{-}[drr] & *{\,}\ar@{-}[dr] &  & *{\,}\ar@{-}[dl] & *{\,}\ar@{-}[dll]\\
 & *{\bullet}\ar@{-}[d] & \ar@{}[l]|{q\,\,\,\,\,\,\,\,\,\,} & \approx &  &  & *{\bullet}\ar@{-}[d] & \ar@{}[l]|{q\circ_{i}p}\\
 & *{\,} &  &  &  &  & *{\,}}
\]

The operad structure of $W(\mathcal{P})$ is given by grafting of
trees, giving the newly arising inner edges length 1. The map $W(\mathcal{P})\rightarrow\mathcal{P}$
is given by setting all lengths to zero (i.e., forget the lengths
and compose in $\mathcal{P}$). 

In \cite{BM W-cons}, it is explained in detail how the above construction
can be performed and studied in the more general context of operads
in any symmetric monoidal category $(\mathcal{E},\otimes,I)$, where
$[0,1]$ is replaced by a suitable ''interval'' $H$ in $\mathcal{E}$.
This is an object $H$ equipped with two ''points'' $\xymatrix{0,1:I\ar@<0.5ex>[r]\ar@<-0.5ex>[r] & H}
$, an augmentation $\epsilon:H\rightarrow I$ satisfying $\epsilon0=id=\epsilon1$,
and a binary operation $\vee:H\otimes H\rightarrow H$ (playing the
role of $max$) which is associative, and for which $0$ is unital
and $1$ is absorbing ($0\vee x=x=x\vee0$ and $1\vee x=1=x\vee1$).
This defines for any operad $\mathcal{P}$ in $\mathcal{E}$ a new
operad $W_{H}(\mathcal{P})$ in $\mathcal{E}$ mapping to $\mathcal{P}$.
The algebras for this operad are up-to-homotopy $\mathcal{P}$-algebras.
For example, one can take for $\mathcal{E}$ the category $Cat$ of
small categories (considered as a model category with weak equivalences
being categorical equivalences, cofibrations those functors that are
injective on objects, and fibrations those functors having the right
lifting property with respect to the functor $0\rightarrow H$, where
$H$ is the groupoid $0\leftrightarrow1$ with two objects and one
isomorphism between them and it also plays the role of the interval). We will examine this possibility below when we consider weak
$n$-categories. 

\begin{example}
Let $[n]$ be the linear tree, viewed as a (discrete) topological
operad. So an $[n]$-algebra consists of a sequence of spaces $X_{0},\cdots X_{n},$
together with maps $f_{ji}:X_{i}\rightarrow X_{j}$ for $i\le j$,
such that $f_{ii}=id$ and \[
\quad\quad\quad\quad\quad\quad\quad\quad\quad\quad f_{kj}\circ f_{ji}=f_{ki}\textrm{ if $i\le j\le k$}\quad\quad\quad\quad\quad\quad\quad\quad\quad\quad\,\,\,\,\,\,\,\,\,\,\,\,\,(1)\]
 A $W([n])$-algebra consists of such a sequence of spaces and maps,
for which (1) holds only up to specified coherent higher homotopies.
Since $W([n])$ is an operad with unary operations only, one can also
think of it as a topological category: it has objects $0,1,\cdots,n$,
and an arrow $i\rightarrow j$ in $W[n]$ is a sequence of ''times''
$t_{i+1},\cdots,t_{j-1}$ (each $t_{k}\in[0,1]$). In other words,
$W[n](i,j)$ is the cube $[0,1]^{j-i-1}$ for $i+1\le j$, a point for $i=j$, and the empty set for $i>j$. Composition is given by juxtaposing two such sequences, putting an
extra time $1$ in the middle: $(t_{i+1},\cdots,t_{j-1}):i\rightarrow j$
and $(t_{j+1},\cdots,t_{k-1}):j\rightarrow k$ compose to give $(t_{i+1},\cdots,t_{j-1},1,t_{j+1},\cdots,t_{k-1})$.
If $\mathcal{C}$ is a category enriched in $Top$ (i.e., a topological
category with discrete set of objects), then the sets of continuous
functors\[
Top(W[n],\mathcal{C})\]
 for varying $n$ define a simplicial set, which is exactly the homotopy
coherent nerve of $\mathcal{C}$, described in \cite{hom lim vogt}.
More generally, if $\mathcal{E}$ is a symmetric monoidal category
with interval $H$, one can construct an $\mathcal{E}$-enriched category
$W_{H}[n]$ with \[
W_{H}[n](i,j)=H^{\otimes_{j-i-1}}\]
 and define for each $\mathcal{E}$-enriched category $\mathcal{C}$
its homotopy coherent nerve $hcN(\mathcal{C})$ as the simplicial
set given by\[
hcN(\mathcal{C})_{n}=\mathcal{E}\textrm{-}Cat(W_{H}[n],\mathcal{C})\]
 that is the set of all $\mathcal{E}$-enriched functors from $W_{H}[n]$
to $\mathcal{C}$. For example, if $\mathcal{E}=Cat$ and $H=0\leftrightarrow1$
as above, then an element of $hcN(\mathcal{C})_{2}$ is given by a
triangle which composes up to a specified invertible 2-cell in $\mathcal{C}$,
\[
\xymatrix{x_{0}\ar[r]\ar[rd] & x_{1}\ar[d]\ar@{}[dl]|<<<<<{\simeq}\\
*{} & x_{2}}
\]

\end{example}
The above generalizes in a completely straightforward way to operads.
Suppose $\mathcal{E}$ and $H$ are as above. Each tree $T$ defines
an operad $\Omega(T)$ in $Set$, which we can view as an operad in
$\mathcal{E}$ (via the functor $Operad\rightarrow Operad(\mathcal{E})$).
Applying the generalized Boardman-Vogt construction yields an operad
$W_{H}(T)$ in $\mathcal{E}$. This construction produces a functor
$\Omega\rightarrow Operad(\mathcal{E})$, which induces an adjunction\[
\xymatrix{Operad(\mathcal{E})\ar@<0.5ex>[r]^{\,\,\,\,\,\, hcN_{d}} & dSet\ar@<0.5ex>[l]^{\,\,\,\,\,\,\,\,\,\,|\cdot|_{H}}}
\]
by Kan extension. For an operad $\mathcal{Q}$ in $\mathcal{E}$ the
dendroidal set $hcN_{d}(\mathcal{Q})$ is called the \emph{homotopy
coherent dendroidal nerve} of $\mathcal{Q}$, and is given explicitly
by \[
hcN_{d}(\mathcal{\mathcal{Q}})_{T}=Hom_{Operad(\mathcal{E})}(W_{H}(T),\mathcal{Q}).\]

\begin{rem}
The functor $|-|_{H}$ is closely related to the $W$- construction
for operads. In fact, if $\mathcal{P}$ is an operad in $Set$, then
the Boardman-Vogt resolution $W_{H}(\mathcal{P}_{\mathcal{E}})$,
of $\mathcal{P}$ viewed as an operad in $\mathcal{E}$, is isomorphic
to the operad $|N_{d}(\mathcal{P})|_{H},$ as follows by direct inspection
of the explicit construction of $W_{H}(\mathcal{P}_{\mathcal{E}})$
in \cite{BM W-cons}. In particular, for an operad $\mathcal{P}$
in $Set$ and an operad $\mathcal{Q}$ in $\mathcal{E}$, there is
a natural bijective correspondence\[
Hom_{Operad(\mathcal{E})}(W_{H}(\mathcal{P}_{\mathcal{E}}),\mathcal{Q})=Hom_{dSet}(N_{d}(\mathcal{P}),hcN_{d}(\mathcal{Q})).\]

\end{rem}
~

\begin{rem}
Consider the special case where $\mathcal{E}$ is the category $Top$
of compactly generated spaces, and $H$ is the unit interval. If $\mathcal{P}$
is a topological operad and $T$ is a tree (an object of $\Omega$),
then the set $hcN_{d}(\mathcal{P})_{T}$ of maps of topological operads
$W_{H}(\Omega(T))\rightarrow\mathcal{P}$ has a natural topology,
as a topological sum of generalized mapping fibrations. (For example,
for the tree $T$ with edges numbered $1,\cdots,5,$
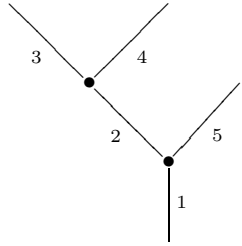
\[
\xymatrix{*{}\ar@{-}[dr]_{3} &  & *{}\ar@{-}[dl]^{4}\\
 & *{\bullet}\ar@{-}[dr]_{2} &  & *{}\ar@{-}[dl]^{5}\\
 &  & *{\bullet}\ar@{-}[d]^{1}\\
 &  & *{}}
\]
$hcN_{d}(\mathcal{P})_{T}$ is the sum, over all $5$-tuples $c_{1},\cdots,c_{5}$
of colours of $\mathcal{P}$, of mapping fibrations of the maps \[
\mathcal{P}(c_{3},c_{4};c_{2})\times\mathcal{P}(c_{2},c_{5};c_{1})\rightarrow\mathcal{P}(c_{3},c_{4},c_{5};c_{1}).)\]
Let $(dTop)^{\delta}$ be the category of dendroidal spaces with discrete
set of vertices. Then $hcN_{d}(\mathcal{P})$ with this topology defines
a functor $hcN_{d}(-):Operad(Top)\rightarrow(dTop)^{\delta}$. This
functor again has a left adjoint $|-|_{H},$which relates to the Boardman-Vogt
resolution of topological operads in the same way as above, by a natural
isomorphism\[
W_{H}(\mathcal{P})\cong|N_{d}(\mathcal{P})|_{H},\]
where $N_{d}$ and $|-|_{H}$ are now viewed as functors between the
categories $Operad(Top)$ and $(dTop)^{\delta}$. 
\end{rem}

\section{The inner Kan condition for dendroidal sets}

Let us begin by recalling some well known facts for simplicial sets
(see for example \cite{calculus of fractions}). Let $\Lambda^{k}[n]\subseteq\Delta[n]$
be the sub-simplicial set of the standard $n$-simplex, defined as
the union of all the faces of $\Delta[n]$ except the one opposite
the $k$-th vertex. A simplicial set $X$ is said to satisfy the Kan
condition, or to be a Kan complex, if for any $n\ge0$ and any $k$
with $0\le k\le n$, any map $\Lambda^{k}[n]\rightarrow X$ can be
extended to a map $\Delta[n]\rightarrow X$. When this is required
for $0<k<n$ only, $X$ is said to be an \emph{inner Kan complex}, or,
to satisfy the \emph{inner Kan condition}. This condition has been
introduced (under the name "restricted Kan condition") by Boardman and Vogt in \cite{BV hom. inv. str.}, while
inner Kan complexes are being studied by Joyal \cite{J Quasi-cat Book}
under the name \emph{quasi-categories}. Observe that the nerve of a category is always an inner Kan complex. In this
section we extend the notion of an inner Kan complex to the context
of dendroidal sets. 

Consider a tree $T$. Recall that a \emph{face} of $T$ is a map $S\rightarrow T$
which corresponds to either contracting an inner edge in $T$ or
pruning an outer cluster in $T$. Those corresponding to an edge
contraction, i.e., $\partial_{e}:T/e\rightarrow T$ for an inner edge
$e$ in $T$, are called \emph{inner faces}. Let $\Lambda^{e}[T]\subseteq\Omega[T]$
be the dendroidal subset of the representable dendroidal set $\Omega[T],$
generated by all the faces of $T$ \emph{except} the inner face $\partial_{e}$.
A dendroidal set $X$ is said to satisfy the \emph{inner Kan condition}
if, for any tree $T$ and any inner edge $e$ in $T$, any map
$\Lambda^{e}[T]\rightarrow X$ extends to a map $\Omega[T]\rightarrow X$
(i.e., to an element in $X_{T}$). A dendroidal set satisfying the
inner Kan condition is also called an \emph{inner Kan complex.}

We now list some examples and properties of dendroidal inner Kan complexes.
Some of the proofs involved are quite technical, and we refer to a
companion paper \cite{MW the second paper} for a detailed exposition
of the proofs. 

\begin{example}

For any operad $\mathcal{P}\in Operad$, (one can easily check
that) the dendroidal nerve $N_{d}(\mathcal{P})$ is an inner Kan complex
(in fact, the extension is unique, and this characterizes those dendroidal
sets that are nerves of operads). 
\end{example}

More generally we have the following.
\begin{prop} 
Let $\mathcal{E}$ be a monoidal model category with a chosen interval $H$. If $\mathcal{P}\in Operad(\mathcal{E})$
is a fibrant operad in the sense that each $\mathcal{P}(c_{1},\cdots,c_{n};c)$
is fibrant, then the homotopy coherent nerve $hcN_{d}(\mathcal{P})$
satisfies the inner Kan condition. 
\end{prop}

A special case of this for simplicial categories was proved in \cite{CorPor}.

\begin{rem}
Inner Kan simplicial sets and inner Kan dendroidal sets are related
as follows. For a simplicial set $X$, the dendroidal set $i_{!}(X)$
is inner Kan iff $X$ is. It is also true that if $Y$ is a dendroidal
set satisfying the inner Kan condition then the simplicial set $i^{*}(Y)$
is again an inner Kan complex. The characterization of the nerves of
operads as those dendroidal sets having unique fillers, is then the
direct analogue of the well known fact that a simplicial set is the
nerve of a category iff it is inner Kan with unique fillers.
\end{rem}

The Grothendieck construction introduced above respects the inner Kan
condition in the following sense.
\begin{prop}
If $X:\mathbb{S}^{op}\rightarrow dSet$
is a diagram of dendroidal sets, each of which is an inner Kan complex,
then the dendroidal set $\int_{\mathbb{S}}X$ is also an inner Kan complex.
\end{prop}

Following \cite{Cisinski}, we call a dendroidal set $X$ \emph{normal}
if, for every object $T$ of $\Omega$ and for every non-degenerate dendrex $x\in X_{T}$, the only
automorphism of $T$ which fixes $x$ is the identify. For example if $X$ is any simplicial set, then $i_{!}(X)$ is normal. And if $\mathcal{P}$ is a $\Sigma$-free operad (i.e., each $\Sigma_{n}$ acts freely), then $N_{d}(\mathcal{P})$ is normal. 

\begin{thm}
Let $K$ be a dendroidal set satisfying the inner Kan condition and
let $X$ be a normal dendroidal set. The dendroidal set $\underline{Hom}_{dSet}(X,K)$
satisfies the inner Kan condition.
\end{thm}
The proof is based on a careful analysis of shuffles of trees, together with the fact that normal dendroidal sets admit a nice skeletal filtration.
This theorem specializes to simplicial sets. Indeed, if $X$ and $K$
are simplicial sets and $K$ is inner Kan, then so is $i_{!}(K)$,
and hence $\underline{Hom}(i_{!}(X),i_{!}(K))$ is a dendroidal inner
Kan complex. Applying $i^{*}$ to it, we see that $5.3(iii)$ implies
that $\underline{Hom}(X,K)$ is a simplicial inner Kan complex. This
simplicial result was already proved by Joyal \cite{J Quasi-cat Book}.
Our proof of Theorem 7.5 thus provides in particular a proof of Joyal's
result, different from the one given in \cite{J Quasi-cat Book} (and similar to the one by \cite{Josh}).

\section{Applications and further developments}

In this last, somewhat speculative section, we would like to point
out some possible further developments of the theory of dendroidal
sets, related to ''weak'' maps between up-to-homotopy algebras,
to enriched and weak higher categories, and to Quillen model categories. 

To begin with, let $\mathcal{P}$ be an operad in $Set$. If $\mathcal{E}$
is a symmetric monoidal model category with a suitable interval $H$,
then $W_{H}(\mathcal{P})$ is an operad in $\mathcal{E}$ whose algebras
are homotopy $\mathcal{P}$-algebras (as mentioned in 6 above). The
maps of $W_{H}(\mathcal{P})$-algebras are maps of homotopy $\mathcal{P}$-algebras
which strictly commute with all higher homotopies, and this is a notion
of map which for many purposes is too restrictive. It is possible
to define a weaker notion of map between homotopy $\mathcal{P}$-algebras,
but then the question arises to what extent these weak maps form a
category. (Boardman and Vogt \cite{BV hom. inv. str.} construct a
''quasi-category'' of weak maps in the context of topological spaces;
in \cite{BM W-cons}, Theorem 6.9, a kind of Segal category of weak
maps is constructed in the context of left proper monoidal model categories;
in \cite{Hess} this question is approached via bimodules). The theory
of dendroidal sets is relevant here. Indeed, $W_{H}(\mathcal{P})$-algebras
in $\mathcal{E}$ are the same thing as operad maps $W_{H}(\mathcal{P})\rightarrow\mathcal{E}$,
or equivalently, as maps of dendroidal sets $N_{d}(\mathcal{P})\rightarrow hcN_{d}(\mathcal{E})$
(see Remark 6.2 above). They thus arise as the vertices of the dendroidal
set \[
\begin{array}{cccc}
   \quad\quad\quad\quad\quad\quad\quad\quad\quad \,\,
 & \underline{Hom}_{dSet}(N_{d}(\mathcal{P}),hcN_{d}(\mathcal{E})). 
 & \quad\quad\quad\quad\quad\quad\quad\quad\quad \,\, &   (1)\end{array}\]
 Dendrices of shape $i[1]$ (where $i:\Delta\rightarrow\Omega$) encode
a suitable notion of weak map, and such weak maps can be composed
(in an up-to-homotopy way) whenever this dendroidal set $(1)$ is
an inner Kan complex. This is the case, for example, when $\mathcal{P}$ is $\Sigma$-free and every object
in $\mathcal{E}$ is fibrant, \emph{c.f.} Proposition 7.2. 

Notice that, more generally, one might consider (weak) $\mathcal{P}$-algebras
with values in any dendroidal set $X$, as vertices of the dendroidal
Hom-set\[
\underline{Hom}_{dSet}(N_{d}(\mathcal{P}),X).\]
If $\mathcal{P}$ is $\Sigma$-free then this dendroidal set is an inner Kan complex whenever $X$ is (Theorem
7.5), in which case maps between $\mathcal{P}$-algebras (again defined
as dendrices of shape $i[1]$) can be composed. The case $X=hcN_{d}(\mathcal{E})$
is the one discussed before. It is also possible to iterate this construction,
and consider for another operad $\mathcal{Q}$ the dendroidal set
\[
\underline{Hom}_{dSet}(N_{d}\mathcal{Q},\underline{Hom}_{dSet}(N_{d}\mathcal{P},X))\]
which is of course isomorphic to\[
\underline{Hom}_{dSet}(N_{d}(\mathcal{P})\otimes N_{d}(\mathcal{Q}),X).\]
This dendroidal set admits a map from \[
\underline{Hom}_{dSet}(N_{d}(\mathcal{P}\otimes_{BV}\mathcal{Q}),X)\]
 but is in general not isomorphic to it, unless $X$ is (the dendroidal
nerve of) an operad. In particular, for the case $X=hcN_{d}(\mathcal{E})$,
one has a map\[
\underline{Hom}(W_{H}(\mathcal{P}\otimes_{BV}\mathcal{Q}),\mathcal{E})\rightarrow\underline{Hom}(|N_{d}(\mathcal{P})\otimes N_{d}(\mathcal{Q})|_{H},\mathcal{E})\]
which gives different but related notions of iterated weak algebras
in $\mathcal{E}$. It would be interesting to compare this to the
work of Dunn, Fiedorowicz, and Vogt on the tensor product of operads
(see e.g., \cite{Dunn,Fiedorowicz}) . (In this context, we should point out
that, up to now, $\mathcal{P}$ and $\mathcal{Q}$ have been operads
in $Set$, but the same applies to topological operads. Indeed, for
the category $Top$ of compactly generated spaces, the homotopy coherent
dendroidal nerve $hcN_{d}(Top)$ with respect to the usual unit interval
is naturally a (large) dendroidal space. If $\mathcal{P}$ is an operad
in $Top$, then homotopy $\mathcal{P}$-algebras in $Top$ are the
vertices of the dendroidal space $\underline{Hom}(N_{d}(\mathcal{P}),hcN_{d}(Top))$,
etc. We expect that (under suitable cofibrancy conditions on $\mathcal{P}$)
this dendroidal space satisfies the inner Kan condition.) 

We would like to consider the special case of the operad $\mathcal{A}_{S}$
whose algebras are categories with a given set $S$ as objects (Example
2.4). Note that this operad is $\Sigma$-free (like any operad obtained by symmetrization, cf. Remark 2.3). For a fixed dendroidal set $X$, one can consider the dendroidal
set \[
\underline{Hom}(N_{d}(\mathcal{A}_{S}),X).\]
By definition, we call its vertices $X$-enriched categories over
$S$. Its dendrices of shape $i[1]$ provide an interpretation of
the notion of ''functor'' between $X$-enriched categories over
$S$. By varying $S$, one obtains a $Set$-indexed diagram of dendroidal
sets, the totalization (see 4.1) of which we denote\[
\begin{array}{cccc}
  \quad\quad\quad\quad\quad\quad\quad\quad\quad\quad\quad\quad\quad\quad &\underline{Cat}(X)   
& \quad\quad\quad\quad\quad\quad\quad\quad\quad\quad\quad\quad\quad\quad & (2)\end{array}\]
By definition, its vertices are categories enriched in $X$, while
its dendrices of shape $i[1]$ are functors between such categories.
In this context, it is relevant to observe that by Theorem 7.5 and
Proposition 7.4, $\underline{Cat}(X)$ is a dendroidal inner Kan complex
whenever $X$ is, so that a composition of functors between $X$-enriched
categories exists. We also note that the construction can be iterated,
so as to form the dendroidal inner Kan complex \[
\underline{Cat}^{2}(X)=\underline{Cat}(\underline{Cat}(X))\]
of $X$-enriched bicategories, and so on. 

Let us consider a few special cases of this construction. First of
all, if $\mathcal{E}$ is a symmetric monoidal category, one can construct
its dendroidal nerve $N_{d}(\mathcal{E})$. The dendroidal set $\underline{Cat}(N_{d}(\mathcal{E}))$
then captures the \emph{usual} notion of $\mathcal{E}$-enriched categories
and functors. More precisely, it is isomorphic to the dendroidal nerve
of the usual monoidal category $Cat(\mathcal{E})$ of $\mathcal{E}$-enriched
categories, \[
\underline{Cat}(N_{d}(\mathcal{E}))\cong N_{d}(Cat(\mathcal{E})),\]
where $Cat(\mathcal{E})$ is considered as an operad via the usual
tensor product of enriched categories. As a particular case, consider
the category $\underline{Cat}$ of small categories with its cartesian
monoidal structure. Then the dendroidal set $Cat^{n}(N_{d}(\underline{Cat})),$
obtained by iterating the construction $n$ times, is the dendroidal
nerve of the category of strict $(n+1)$-categories. It also encodes all
higher structure of functors, natural transformations, modifications,
and so on.

If $\mathcal{E}$ is a monoidal model category with a suitable interval
$H$, one can consider categories enriched in the homotopy coherent
nerve $hcN_{d}(\mathcal{E})$ (defined in terms of $H$). For example,
if $\mathcal{E}$ is the category of chain complexes over a ring $R$
(with the projective model structure and the usual interval $H$ of
normalized chains on the standard $1$-simplex), then $\underline{Cat}(hcN_{d}(\mathcal{E}))$
is a dendroidal inner Kan complex whose vertices are precisely $A_{\infty}$-categories
(\cite{Fukaya,Lefevre,Lybashenko}). As another example, let $\mathcal{E}=Top$
with the unit interval, and consider for the one-point set $*$ the
dendroidal inner Kan complex\[
A_{\infty}=\underline{Hom}(N_{d}(\mathcal{A}_{*}),hcN_{d}(\mathcal{E})).\]
The vertices of this dendroidal set are precisely $A_{\infty}$-spaces,
while dendrices of more general shapes encode operations between $A_{\infty}$-spaces.
Again, the construction can be iterated to form dendroidal inner Kan
complexes $A_{\infty}^{(1)}=A_{\infty}$ and \[
A_{\infty}^{(n+1)}=\underline{Hom}(N_{d}(\mathcal{A}_{*}),A_{\infty}^{(n)}).\]
It would be interesting to study the relation between $A_{\infty}^{(n)}$
and $n$-fold loop spaces in topology \cite{Dunn, iterated loops}. 

Finally, the category $\underline{Cat}$ of small categories itself
is a monoidal model category with interval $H$ as in 6 above, and
\[
\underline{Hom}(N_{d}(\mathcal{A}_{S}),hcN_{d}(\underline{Cat}))\]
is a dendroidal inner Kan complex capturing the notion of a \emph{bicategory}
with $S$ as set of objects \cite{Benabou} (or rather the notion of an unbiased bicategory \cite{Leinster book}). This construction can
again be iterated. For example, the dendroidal inner Kan complex $\underline{Hom}(N_{d}(\mathcal{A}_{*}),\underline{Hom}(N_{d}(\mathcal{A}_{*}),hcN_{d}(\underline{Cat})))$
captures braided monoidal categories (and all higher maps between
them). The above construction of categories enriched in $\mathcal{E}$
yields, by considering $\mathcal{E}=\underline{Cat}$, an inductive
definition of weak $n$-categories. More precisely, let $WCat_{1}=\underline{Cat}$
and $ $for $n>1$ let \[
WCat_{n}=\underline{Cat}^{n-1}(hcN_{d}(\underline{Cat})).\]
For each $n\ge1$, $WCat_{n}$ is a dendroidal inner Kan complex. Its
vertices are weak $n$-categories of a special kind. (They have an
underlying strict category of $1$-cells, and for any two objects
$x$ and $y$, the same is true at level $n-1$ for the dendroidal
set $Hom(x,y)$). There are many alternative notions of weak $n$-categories
in the literature (see \cite{survey of n-cat} for a survey of 10
such definitions and \cite{Baez Weak n-cat} for a more general discussion
of weak $n$-categories), and we expect that for
any reasonable notion, a weak $n$-category can be ''strictified''
to a weak $n$-category in our sense. 

Finally, we would like to say a few words about possible Quillen model
structures on dendroidal Sets. Recall from \cite{J Quasi-cat Book,Lurie}
that there is a Quillen model structure on simplicial sets, in which
the inner Kan complexes are exactly the fibrant objects. This model
structure is related to the ''folk'' monoidal model structure on
$\underline{Cat}$ already mentioned above, in which the weak equivalences
are the equivalences of categories and the cofibrations are the functors
which are injective on objects. Indeed, a map $X\rightarrow Y$ between
simplicial sets is a weak equivalence in Joyal's model structure iff
for every simplicial inner Kan complex $K$, the map $\tau(K^{Y})\rightarrow\tau(K^{X})$
is an equivalence of categories (here $\tau:sSet\rightarrow Cat$
is the functor discussed in 4). The analog of Theorem 7.5 for simplicial
sets, which states that $K^{X}$and $K^{Y}$ are again inner Kan complexes,
plays an important role in Joyal's model structure. 

The folk model structure on $\underline{Cat}$ generalizes without
much effort to one on (coloured) operads, in which a map $f:\mathcal{Q}\rightarrow\mathcal{P}$,
from an operad $\mathcal{Q}$ on $D$ to an operad $\mathcal{P}$
on $C$ (as in Section 2) is a weak equivalence iff $j^{*}(f):j^{*}\mathcal{Q}\rightarrow j^{*}\mathcal{P}$
is an equivalence of categories, and moreover $f$ induces a bijection
\[
\mathcal{Q}(d_{1},\cdots,d_{n};d)\rightarrow\mathcal{P}(fd_{1},\cdots,fd_{n};fd)\]
 for any sequence $d_{1},\cdots,d_{n},d$ of colours in $D$. We conjecture
that the inner Kan complexes are the fibrant objects in a model structure
on dendroidal sets in which a map $X\rightarrow Y$ is a weak equivalence
iff, for any dendroidal inner Kan complex $K$, the map\[
\tau_{d}(\underline{Hom}_{dSet}(Y,K)\rightarrow\tau_{d}(\underline{Hom}_{dSet}(X,K))\]
is a weak equivalence of operads. Theorem 7.5 should be a substantial
step towards a proof of this conjecture.

\section{Appendix: The tensor product of dendroidal objects }

Let $\mathcal{E}$ be a symmetric monoidal category. The category
of dendroidal objects in $\mathcal{E}$ is the functor category $\mathcal{E}^{\Omega^{op}}$,
which we denote by $d\mathcal{E}$. This category has a Boardman-Vogt
style tensor product, and a corresponding internal $\underline{Hom}$
whenever $\mathcal{E}$ itself is closed. The construction and its
basic properties are explained most easily after recalling some basic
facts about ''enriched Kan extensions'', so we'll do that first.
None of the material in 9.1 is really new, and we refer the reader
to \cite{Kelly ECT} for more background.

\subsection{Enriched Kan extensions}

We begin by developing a bit of formalism similar to the language
of rings and bimodules. Let $\mathcal{E}$ be a symmetric monoidal
category, and let $\mathcal{S}$ be any $\mathcal{E}$-enriched category.
Suppose $\mathcal{S}$ is \emph{tensored} over $\mathcal{E}$. (This
means that one can construct an object $E\otimes S$ in $\mathcal{S}$
for $E$ in $\mathcal{E}$ and $S$ in $\mathcal{S}$, with the property
that there is a natural $\underline{Hom}$-tensor correspondence between
maps $E\otimes S\rightarrow T$ in $\mathcal{S}$ and $E\rightarrow\underline{Hom}(S,T)$
in $\mathcal{E}$; see \cite{Kelly ECT} for a formal definition.
For small categories $\mathbb{A}$ and $\mathbb{B}$ (in $Set$),
we write\[
_{\mathbb{A}}\mathcal{E}_{\mathbb{B}}=\mathcal{E}^{\mathbb{B}^{op}\times\mathbb{A}}\]
 for the category of functors $\mathbb{B}^{op}\times\mathbb{A}\rightarrow\mathcal{E}$.
For objects $X\in$$_{\mathbb{A}}\mathcal{E}_{\mathbb{B}}$ and $A$
in $\mathbb{A}$, $B$ in $\mathbb{B}$, we write\[
_{A}X_{B}=X(B,A)\in\mathcal{E}\]
 for the value at $(B,A)$. Also, if $\mathbb{A}$ or $\mathbb{B}$
is the trivial category $\star$ we delete it from the notation. So\[
_{\mathbb{A}}\mathcal{E}_{\star}=\,_{\mathbb{A}}\mathcal{E}=\mathcal{E}^{\mathbb{A}},\,\,\,_{\star}\mathcal{E}_{\mathbb{B}}=\mathcal{E}_{\mathbb{B}}=\mathcal{E}^{\mathbb{B}^{op}}.\]
 Now assume $\mathcal{E}$ has small limits and $\mathcal{S}$ has
small colimits. There is a tensor product functor \[
\begin{array}{cccc}
\quad\quad\quad\quad\quad\quad\quad\quad\quad\quad& \otimes_{\mathbb{B}}=\,_{\mathbb{C}}\mathcal{E}_{\mathbb{B}}\times\,_{\mathbb{B}}\mathcal{S}_{\mathbb{A}}\rightarrow\,_{\mathbb{C}}\mathcal{S}_{\mathbb{A}}  & 
\quad\quad\quad\quad\quad\quad\quad\quad\quad\quad& (1)\end{array}\]
 defined for $E$ in $_{\mathbb{C}}\mathcal{E}_{\mathbb{B}}$ and
$S$ in $_{\mathbb{B}}\mathcal{S}_{\mathbb{A}}$, by the usual coequalizer\[
\xymatrix{_{C}(E\otimes_{\mathbb{B}}\mathcal{S})_{A} & \coprod_{B}(_{C}E_{B})\otimes(_{B}S_{A})\ar[l] & \coprod_{B\rightarrow B'}(_{C}E_{B'})\otimes(_{B}S_{A})\ar@<+2pt>[l]\ar@<-2pt>[l]}
\]
 for any two objects $C\in\mathbb{C}$, $A\in\mathbb{A}$. This tensor
product has a corresponding internal $\underline{Hom}$, \[
\begin{array}{cccc}
\quad\quad\quad\quad\quad\quad\quad\quad\quad\,\,& \underline{Hom}_{\mathbb{A}}:\,_{\mathbb{B}}\mathcal{S}_{\mathbb{A}}\times\,_{\mathbb{C}}\mathcal{S}_{\mathbb{A}}\rightarrow\,_{\mathbb{C}}\mathcal{E}_{\mathbb{B}},   &
\quad\quad\quad\quad\quad\quad\quad\quad\quad\,\,& (2)\end{array}\]
 satisfying the usual adjunction property stating a bijective correspondence
between maps\[
\begin{array}{cc}
E\otimes_{\mathbb{B}}S\rightarrow T & \textrm{in $_{\mathbb{C}}\mathcal{S}_{\mathbb{A}}$}\end{array}\]
and maps\[
\begin{array}{cc}
E\rightarrow\underline{Hom}_{\mathbb{A}}(S,T) & \textrm{in $_{\mathbb{C}}\mathcal{E}_{\mathbb{B}}$}.\end{array}\]
 We point out two special cases of this Hom-tensor correspondence.
First, if $F$ is an element of $_{\mathbb{B}}\mathcal{S}_{\mathbb{A}}$,
i.e., $F:\mathbb{B}\rightarrow\mathcal{S}_{\mathbb{A}},$ then we
obtain adjoint functors\[
\xymatrix{\quad\quad\quad\quad\quad\quad\quad\quad\quad\quad\quad\quad\,\,\, f_{!}:\mathcal{E}_{\mathbb{B}}\ar@<2pt>[r] & \mathcal{S}_{\mathbb{A}}:f^{*}\ar@<2pt>[l] 
&\quad\quad\quad\quad\quad\quad\quad\quad\quad\,\,\, (3)}
\]
 defined in terms of the previous functors for the special case $\star=\mathbb{C}$,
by\[
f_{!}(E)=E\otimes_{\mathbb{B}}F\,\,\,\,\,\, f^{*}(S)=\underline{Hom}_{\mathbb{A}}(F,S)=\underline{Hom}(F,S).\]
 These functors $f^{*}$ and $f_{!}$ are the right and left Kan extensions
along $F$. Secondly, there are ''external'' tensor and Hom functors\[
\begin{array}{cccc}
\quad\quad\quad\quad\quad\quad\quad\quad\quad\quad\quad&\underline{\otimes}:\mathcal{E}_{\mathbb{C}}\times\mathcal{S}_{\mathbb{A}}\rightarrow\mathcal{S}_{\mathbb{C}\times\mathbb{A}} &    
\quad\quad\quad\quad\quad\quad\quad\quad\quad\quad\quad& (4)\end{array}\]
 \[
\begin{array}{cccc}
\quad\quad\quad\quad\quad\quad\quad\quad\quad\quad&\underline{Hom}_{\mathbb{A}}:\mathcal{S}_{\mathbb{A}}\times\mathcal{S}_{\mathbb{C}\times\mathbb{A}}\rightarrow\mathcal{E}_{\mathbb{C}}  &  
\quad\quad\quad\quad\quad\quad\quad\quad\quad\quad& (5)\end{array}\]
 for which there is a natural correspondence between maps \[
\begin{array}{cc}
X\underline{\otimes}Y\rightarrow Z & \textrm{in }\mathcal{S}_{\mathbb{C}\times\mathbb{A}}\end{array}\]
 and\[
\begin{array}{cc}
X\rightarrow\underline{Hom}_{\mathbb{A}}(Y,Z) & \textrm{in }\mathcal{E}_{\mathbb{C}}\end{array}\]
 Indeed, this is the special case where $\mathbb{B}=\star$ while
$\mathbb{C}$ is replaces by $\mathbb{C}^{op}$, so that (1) and (2)
can be rewritten as $\otimes:\,_{\mathbb{C}^{op}}\mathcal{E}\times\mathcal{E}_{\mathbb{A}}\rightarrow\,_{\mathbb{C}^{op}}\mathcal{E}_{\mathbb{A}}$
and $\underline{Hom}_{\mathbb{A}}:\mathcal{S}_{\mathbb{A}}\times\,_{\mathbb{C}^{op}}\mathcal{S}_{\mathbb{A}}\rightarrow\,_{\mathbb{C}^{op}}\mathcal{E}$,
defining (4) and (5).

Now consider a functor $F:\mathbb{A}\times\mathbb{A}\rightarrow\mathcal{S}_{\mathbb{A}}$,
i.e., $F\in_{\mathbb{A}\times\mathbb{A}}\mathcal{S}_{\mathbb{A}}$.
Then by Kan extension we have a functor \[
\xymatrix{\mathcal{E}_{\mathbb{A}}\times\mathcal{S}_{\mathbb{A}}\ar[r]^{\underline{\otimes}} & \mathcal{S}_{\mathbb{A}\times\mathbb{A}}\ar[r]^{f_{!}} & \mathcal{S}_{\mathbb{A}}}
\]
 which we write as $\otimes^{(F)}$; so\[
E\otimes^{(F)}S=f_{!}(E\underline{\otimes}S)=(E\underline{\otimes}S)\otimes_{\mathbb{A}\times\mathbb{A}}F.\]
The above discussion also yields a corresponding Hom-functor, denoted\[
\underline{Hom}^{(F)}:\mathcal{S}_{\mathbb{A}}\times\mathcal{S}_{\mathbb{A}}\rightarrow\mathcal{E}_{\mathbb{A}},\]
 for which there is a bijective correspondence between maps\[
E\otimes^{(F)}S\rightarrow T\textrm{\,\,\,\,\,\,(in $\mathcal{S}_{\mathbb{A}}$)}\]
and maps\[
E\rightarrow\underline{Hom}^{(F)}(S,T)\textrm{\,\,\,\,\,(in $\mathcal{E}_{\mathbb{A}}$).}\]
 Indeed, one can simply define $\underline{Hom}^{(F)}$ in terms of
the earlier $\underline{Hom}$ and the right adjoint $f^{*}$, as
\[
\underline{Hom}^{(F)}(S,T)=\underline{Hom}_{\mathcal{A}}(S,f^{*}T)\]

\subsection{Monoidal closed structure of $d\mathcal{E}$}

Let us now consider a complete and cocomplete symmetric closed monoidal
category $\mathcal{E}$, and the category $d\mathcal{E}$ of dendroidal
object in $\mathcal{E}$. Let $\mathbb{A}=\Omega$, let $\mathcal{S}=\mathcal{E}$,
and let $F=BV$ be the Boardman-Vogt tensor product of Hopf operads
in $\mathcal{E}$, restricted to operads coming from $\Omega$:\[
\xymatrix{BV:\Omega\times\Omega\ar[r] & Operad(\mathcal{E})\ar[r]^{\,\,\,\,\,\,\,\,\,\,\, N_{d}} & d\mathcal{E}.}
\]
 Then the last construction of 9.1 yields a functor\[
\otimes^{(BV)}:d\mathcal{E}\times d\mathcal{E}\rightarrow d\mathcal{E}\]
 and a corresponding Hom-functor\[
\underline{Hom}^{(BV)}:d\mathcal{E}\times d\mathcal{E}\rightarrow d\mathcal{E}\]
 satisfying the usual properties, and \emph{making} $d\mathcal{E}$
\emph{into a closed symmetric monoidal category. }

\end{document}